\documentclass[10pt]{amsart}

\usepackage[latin1]{inputenc}
\usepackage{amsmath}
\usepackage{amsfonts}
\usepackage{amssymb}
\usepackage{ytableau}
\usepackage[mathscr]{eucal}
\usepackage{diagrams}
\usepackage{xy,color,graphicx}
\usepackage{hyperref}
\xyoption{all}
\usepackage{array}
\usepackage{enumerate}
\usepackage{url}

\newcommand{\C}{\mathbb{C}}

\newcommand{\Af}{\mathbb{A}}

\newtheorem{theorem}{Theorem}

\title{virtual Representation motives}
\author{Lieven Le Bruyn} 
\address{Department Mathematics, University  of Antwerp ,
 Middelheimlaan 1, B-2020 Antwerp (Belgium) {\tt lieven.lebruyn@uantwerpen.be}}

\begin{document}
\sloppy

\maketitle

\begin{abstract}
Principal $GL_n$-bundles (aka vector bundles) are locally trivial in the Zariski topology, whereas principal $PGL_n$-bundles (aka Azumaya algebras) are not, to the delight of every non-commutative algebraist. Still, this makes the calculation of motives of representation schemes of algebras next to impossible. In very special cases, Brauer-Severi schemes (and their motives) can be used to tackle this problem inductively. We illustrate this in the case of certain superpotential algebras.\end{abstract}

\vskip 4mm

The scheme of $n$-dimensional representations $\mathbf{rep}_n R$ of a finitely presented non-commutative algebra $R$
\[
R=\frac{\C \langle X_1,\hdots,X_m \rangle}{(F_1,\hdots,F_k)} \]
is by definition the zero-set in $\Af^{mn^2} = M_n(\C)^{\oplus m}$ of $kn^2$ equations
\[
F_l(A_1,\hdots,A_m)_{ij} = 0 \qquad 1 \leq l \leq k,~1 \leq i,j \leq n \]
where $A_u$ is the generic $n \times n$ matrix, that is, the $n \times n$ matrix with entries the coordinates in $\Af^{mn^2} = M_n(\C)^{\oplus m}$ corresponding to the $u$-th factor. Very little is known about the global structure of these representation schemes. Here we will be interested in their (virtual) motives.

\section{Motives}

Motives are best thought of as a Lego-version of varieties. That is, for every reduced $\C$-variety $X$ we have, up to Zariski isomorphism, one block $[ X ]$, called the {\em motive} of $X$, which are elements in the ring of (naive) motives $\mathbf{Mot}_{\C}$, in which addition and multiplication are subject to the following rules:
\begin{itemize}
\item{if $Y$ is a Zariski closed subvariety of $X$, then we have a 'scissor-relation'
\[
[X] = [X-Y]+[Y] \]
allowing us to slice up a variety is locally closed parts and compute its motive by adding up these smaller blocks.}
\item{if $X$ is a fiber bundle in the Zariski topology with base $Y$ and fiber $F$, then we have a factoring-relation
\[
[X] = [Y] \times [ F] \]
giving us in particular that the motive of a product is the product of the motives.}
\end{itemize}
The {\em Lefschetz-motive} $\mathbb{L} = [ \Af^1_{\C} ]$ is the motive of the affine line, and all varieties allowing a cell-decomposition (such as projective spaces or Grassmannians) have therefore as their motive a polynomial in $\mathbb{L}$. It is quite easy to verify that
\[
 [ \mathbb{P}^n ] = \frac{\mathbb{L}^{n+1}-1}{\mathbb{L}-1}, \quad [ GL_n ] = \prod_{k=0}^{n-1} (\mathbb{L}^n - \mathbb{L}^k), \quad [ Gr(k,n) ] = \frac{[ GL_n]}{[ GL_k ] [ GL_{n-k}] \mathbb{L}^{k(n-k)}}  \]

Motives can be calculated either via geometric insight or by laborious algebraic manipulations. Consider for example a smooth affine quadric in $\Af^3$
\[
Q = \mathbb{V}(xy-z^2-1) \subseteq \Af^3 \]
To a geometer  this is the affine piece of a smooth quadric in $\mathbb{P}^3_{\C}$ which she knows to be isomorphic to $\mathbb{P}^1 \times \mathbb{P}^1$ from which she has to remove the intersection with the hyperplane at infinity, which is a $\mathbb{P}^1$, so its motive must be the difference
\[
[ Q ] = (\mathbb{L}+1) (\mathbb{L}+1) - (\mathbb{L}+1) = \mathbb{L}^2 + \mathbb{L} \]

An algebraist would chop up the variety in smaller pieces by localization and eliminating variables. On the open piece $x \not= 0$ he can eliminate
\[
y = \frac{z^2+1}{x} \quad \text{so} \quad [x \not= 0 ] = \mathbb{L} (\mathbb{L}-1) \]
because $z$ is a free variable and $x$ cannot be zero. On the complement $x=0$ the equation becomes $z^2=1$ and therefore its motive is
\[
[ x = 0 ] = 2 \mathbb{L} \]
because $z = \pm 1$ whereas $y$ is a free variable. Adding these contribution he gets
\[
[Q] = [ x \not= 0]+[x=0] = \mathbb{L}(\mathbb{L}-1) + 2 \mathbb{L} = \mathbb{L}^2+\mathbb{L} \]

\vskip 3mm 

Similarly, for a singular affine quadric in $\C^3$
\[
C = \mathbb{V}(xy-z^2) \subseteq \C^3 \]
our geometer will view this as a cone over a smooth conic in $\mathbb{P}^2$ so would think of it as $\mathbb{P}^1 \times \C^* \sqcup \{ top \}$ with corresponding motive
\[
[C] = (\mathbb{L}+1)(\mathbb{L}-1) + 1 = \mathbb{L}^2 \]
The algebraist would again decompose into
\[
[ x \not= 0] = [ y = \frac{z}{2}, x \not= 0] = \mathbb{L}(\mathbb{L}-1) \]
and
\[
[ x = 0 ] = [z = 0, x=0] = \mathbb{L}  \]
because $y$ is still a free variable, giving the same answer $\mathbb{L}^2$.

\section{$\neg$~Luna}

Let us return to the representation motives $[ \mathbf{rep}_n R ]$.  An evident approach to representation varieties is via invariant theory. That is, consider the action of $GL_n$ on $\mathbf{rep}_n R$ via basechange (conjugation), and study the corresponding quotient-variety 
\[
\pi : \mathbf{rep}_n R \rOnto \mathbf{rep}_n R / GL_n = \mathbf{iss}_n. R \]
By Mumford's GIT the points of $\mathbf{iss}_n R$ correspond to closed orbits which by a result of Michael Artin we know are the isomorphism classes of $n$-dimensional {\em semi-simple} representations of $R$. The quotient map $\pi$ sends an $n$-dimensional $R$-module to the isomorphism class of the direct sum of its Jordan-H\"older components.

\vskip 3mm

So, we might try to decompose the quotient variety $\mathbf{iss}_n R$ according to different representation types $\tau$ of semi-simples and calculate the motives $[\mathbf{iss}_n R (\tau) ]$ of these so called {\em Luna strata}. 

In rare cases, for example if the representation scheme $\mathbf{rep}_n R$ is a smooth variety, one can show that the fibers $\pi^{-1}(S_{\tau})$ are isomorphic for all $S_{\tau} \in \mathbf{iss}_n R (\tau)$, so we might hope in such cases to arrive at the representation motive via
\[
[ \mathbf{rep}_n R ] \overset{?}{=} \sum_{\tau} [\pi^{-1}(S_{\tau})] [ \mathbf{iss}_n R (\tau) ] \]
Even in the simplest of cases one obtains nonsense.

\vskip 3mm

Let $R = \C [x]$, then clearly, $\mathbf{rep}_2 \C[x] = M_2(\C)$ determined by the matrix-image of $x$, and the $GL_2$-action is by conjugation. The quotient map assigns to a matrix the coefficients of its characteristic polynomial, so the quotient map is
\[
\mathbf{rep}_2 \C[x] = M_2(\C) \rOnto^{\pi} \C^2 = \mathbf{iss}_2 \C[x] \qquad A \mapsto (Tr(A),Det(A)) \]

Semi-simple matrices are the diagonalisable ones, so there are two representations types of semi-simples: $\tau_1$ with two distinct eigenvalues and $\tau_2$ with two equal eigenvalues. The second stratum is determined by the closed subvariety $\mathbb{V}(Det-Tr^2)$ which is a smooth conic in $\C^2$. Therefore,
\[
[ \mathbf{iss}_2 \C[x](\tau_2) ] = \mathbb{L} \quad \text{and} \quad [ \mathbf{iss}_2 \C[x](\tau_1) ] = \mathbb{L}^2 - \mathbb{L} \]

The fibers $\pi^{-1}(S_{\tau_1})$  are the closed orbits 
\[
\mathcal{O}(\begin{bmatrix} \lambda_1 & 0 \\ 0 & \lambda_2 \end{bmatrix}) \]
 which are all a smooth affine quadric in $\C^3$ and therefore $[ \pi^{-1}(S_{\tau_1} )] = \mathbb{L}^2+\mathbb{L}$.
 
 On the other hand, the fibers $\pi^{-1}(S_{\tau_2})$ are the orbit-closures
 \[
 \overline{\mathcal{O}(\begin{bmatrix} \lambda & 1 \\ 0 & \lambda \end{bmatrix})} \]
 which are singular affine quadrics in $\C^3$, with the top corresponding to the diagonal matrix. Therefore, $[ \pi^{-1}(S_{\tau_2}] = \mathbb{L}^2$. The Luna stratification approach gives us
 \[
 \sum_{\tau} [ \pi^{-1}(S)_{\tau}][\mathbf{iss}_2 \C[x] (\tau)] = (\mathbb{L}^2+\mathbb{L})(\mathbb{L}^2-\mathbb{L}) + \mathbb{L}^2.\mathbb{L} = \mathbb{L}^4+\mathbb{L}^3-\mathbb{L}^2 \]
 which is clearly different from $[\mathbf{rep}_2 \C[x]] = \mathbb{L}^4$.

\vskip 3mm

What went wrong here is that the fibrations considered are fibrations locally trivial in the {\em \'etale topology} but not in the Zariski topology.
Going from coefficients of the characteristic polynomial to eigenvalues involves taking roots, which are typical examples of \'etale extensions, but of course not isomorphisms. 

 And we can't allow \'etale isomorphisms in defining motives because then the ring of motives would become the trivial ring.

\vskip 3mm

Another way to explain this difficulty is to observe that the $GL_n$-action on $\mathbf{rep}_n R$ is actually an action of $PGL_n$ and that there is a huge difference between these two groups when it comes to fibrations. $GL_n$ is a {\em special} group meaning that all \'etale  principal fibrations are in fact Zariski fibrations. Principal $GL_n$-fibrations over an affine scheme $X$ correspond to rank $n$ projective modules over $\C[X]$.

On the other hand, principal \'etale $PGL_n$-fibrations over $X$ correspond to {\em Azumaya algebras} $A$ over $\C[X]$, that is, algebras $A$ which are projective modules of rank $n^2$ over their center $\C[X]$ such that
\[
A \otimes_{\C[X]} A^{op} \simeq End_{\C[X]}(A) \]
The correspondence is given by assigning to an Azumaya algebra $A$ over $X$ its representation scheme. Then, the quotient map
\[
\mathbf{rep}_n A \rOnto^{\pi} X = \mathbf{iss}_n A \]
is the corresponding principal $PGL_n$-fibration.

The principal Zariski $PGL_n$-fibrations correspond to the {\em trivial} Azumaya algebras, that is, those of the form $A = End_{\C[X]}(P)$ where $P$ is a projective $\C[X]$-module of rank $n$.

This distinction between \'etale and Zariski principal $PGL_n$-fibrations is at the very heart of non-commutative algebra. The obstruction to all Azumaya algebras over $X$ being trivial is measured by an important invariant, the {\em Brauer group} $Br(X)$ of $X$.

\section{Framing}

We have seen that we cannot use the Luna stratification approach in order to compute representation motives, caused by the fact that the acting group on representation schemes is $PGL_n$ rather than $GL_n$.

To bypass this problem we might try to replace the action of $PGL_n$ by one of $GL_n$.
 One way to achieve this is by a process called {\em framing}.
 
 \vskip 3mm

Instead of $\mathbf{rep}_n R$ we consider the product $\mathbf{rep}_n R \times \C^n$ and the action of $GL_n$ on it defined by
\[
g.(\phi,v) = (g.\phi.g^{-1},g.v) \]
As long as $v \not= 0$ we see that non-trivial central elements act non-trivially on the second factor, so this is a genuine $GL_n$-action. 

So, we might try a stratification strategy on the $GL_n$-variety $\mathbf{rep}_n R \times (\C^n - \{ 0 \})$ in order to compute its motive, which is $(\mathbb{L}^n-1)[ \mathbf{rep}_n R]$, by summing over the different strata.

\vskip 3mm

Let us first consider the case when  $A$ is an Azumaya algebra. Then, $GL_n$ acts freely on $\mathbf{rep}_n A \times (\C^n- \{ 0 \})$, and so the corresponding quotient map is a principal $GL_n$-fibration
\[
\pi : \mathbf{rep}_n A \times \C^n- \{ 0 \} \rOnto (\mathbf{rep}_n A \times \C^n- \{ 0 \})/GL_n = \mathbf{BS}_n(A) \]
where $\mathbf{BS}_n(A)$ is called the {\em Brauer-Severi variety} of the Azumaya algebra $A$. As this time the quotient map is a Zariski fibration we have
\[
[ \mathbf{rep}_n(A) ] (\mathbb{L}^n -1) = [\mathbf{BS}_n(A)] [GL_n] \]
That is, we can compute the representation motive of $A$ if we can compute its Brauer-Severi motive.

In the trivial case when $A=M_n(\C)$ we have that $\mathbf{rep}_n M_n(\C) = PGL_n$ and $\mathbf{BS}_n(M_n(\C))=\mathbb{P}^{n-1}$ so the above equality reduces to
\[
[ PGL_n ] (\mathbb{L}^n-1) =  [ \mathbb{P}^{n-1}] [GL_n] \quad \text{that is} \quad [\mathbb{P}^{n-1}] = \frac{\mathbb{L}^n-1}{\mathbb{L}-1} \]

\vskip 3mm

For arbitrary algebras $R$ the situation is of course more complicated, but we can use the above idea to compute representation motives inductively from knowledge of motives of {\em generalised} Brauer-Severi varieties. 

\vskip 3mm

In the product $\mathbf{rep}_n R \times (\C^n - \{ 0 \})$ let us consider the Zariski open subset of {\em stable couples}
\[
\mathbf{S}_{n,n}(R) = \{ (\phi,v)~|~\phi(R)v = \C^n \} \]
on which $GL_n$ acts freely, so we have a principal $GL_n$-fibration
\[
\mathbf{S}_{n,n}(R) \rOnto \mathbf{S}_{n,n}(R)/GL_n = \mathbf{BS}_n(R) \]
with $\mathbf{BS}_n(R)$ the {\em $n$-the Brauer-Severi variety of $R$} as introduced by Michel Van den Bergh\footnote{M. Van den Bergh, {\em The Brauer-Severi scheme of the trace ring of generic matrices}, NATO ASI Vol. 233, 333-338 (1987)}.

\vskip 3mm

We can decompose the product $\mathbf{rep}_n R \times ( \C^n - \{ 0 \})$ into the locally closed strata
\[
\mathbf{S}_{n,k}(R) = \{ (\phi,v)~|~dim_{\C} \phi(R)v = k \} \]
giving us this motivic equality
\[
(\mathbb{L}^n-1) [ \mathbf{rep}_n R ] = \sum_{k=1}^n [ \mathbf{S}_{n,k}(R) ] \quad \text{with} \quad [ \mathbf{S}_{n,n}(R) ] = [ \mathbf{BS}_n(R) ] [ GL_n] \]

\vskip 3mm

In order to calculate the motives of the intermediate srata $\mathbf{S}_{n,k}(R)$ with $1 \leq k \leq n-1$ consider the map $\psi$ sending a couple $(\phi,v)$ to the $k$-dimensional subspace $V=\phi(R).v$ of $\C^n$
\[
\psi : \mathbf{S}_{n,k} \rOnto Gr(k,n) \]
To compute the fiber $\psi^{-1}(V)$ take a basis of $V$ and extend this to a basis for $\C^n$, then with respect to this basis, any couple in the fiber can be written as
\[
(\phi,v) = (\begin{bmatrix} \phi_1 & e \\ 0 & \phi_2 \end{bmatrix},\begin{bmatrix} w \\ 0 \end{bmatrix}) \]
with $(\phi_1,w) \in \mathbf{S}_{k,k}(R)$ and $\phi_2 \in \mathbf{rep}_{n-k}(R)$ and $e \in Ext^1_R(\phi_2,\phi_1)$ an extension of the two representations.

\vskip 3mm

In extremely rare situations it may happen that this extension-space is of constant dimension, say $d$, along $\mathbf{S}_{n,k}(R)$, which would then allow us to compute
\begin{eqnarray*}
[ \mathbf{S}_{n,k}(R) ]  = &  \mathbb{L}^d [ Gr(k,n) ] [ \mathbf{S}_{k,k}(R) ] [ \mathbf{rep}_{n-k}(R)] \\  = &  \mathbb{L}^d [ Gr(k,n) ] [ GL_k ] [\mathbf{BS}_k(R)]  [ \mathbf{rep}_{n-k}(R)] \end{eqnarray*}
If we were so lucky for this to hold for all intermediate strata, we would then be able to compute the representation motive $[\mathbf{rep}_n(R)]$ inductively from knowledge of the representation motives $[ \mathbf{rep}_k(R)]$ for $k < n$ and the Brauer-Severi motives $[ \mathbf{BS}_k(R)]$ for $k \leq n$.

\vskip 3mm

Clearly, one would expect this extension condition to hold only for algebras close to free- or quiver-algebras, and not in more interesting situations. Surprisingly, one can reduce to the almost free setting in the case of {\em superpotential algebras}.

\section{Superpotentials}

A {\em superpotential} is a non-commutative homogeneous word $W \in \C \langle X_1,\hdots,X_m \rangle_d$ of degree $d$ in $m$ variables. It determines a Chern-Simons functional
\[
Tr(W) : M_n(\C)^{\oplus m} \rTo \C \qquad (M_1,\hdots,M_n) = Tr(W(M_1,\hdots,M_n)) \]
For ringtheorists the relevant fact is that the degeneracy locus of this map
\[
\{ d Tr(W)=0 \} = \mathbf{rep}_n R_W \quad \text{where} \quad R_W = \frac{\C \langle X_1,\hdots,X_m \rangle}{(\partial_{X_1},\hdots,\partial_{X_m})} \]
is the representation variety of the corresponding Jacobi algebra where the $\partial_{X_i}$ are the cyclic derivatives with respect to $X_i$.

\vskip 3mm

As an example, take $W=aXYZ+bXZY+\frac{c}{3}(X^3+Y^3+Z^3)$, then these cyclic derivatives are
\[
\begin{cases} \partial_X~:~aYZ+bZY+cX^2 \\
\partial_Y~:~aZX+bXZ+cY^2 \\
\partial_Z~:~aXY+bYX+cZ^2 \\
\end{cases}
\]
so the Jacobian algebra $R_W$ is the $3$-dimensional Sklyanin algebra.

\vskip 3mm
The fibers of the Chern-Simons functional $M_n(\lambda) = Tr(W)^{-1}(\lambda)$ for all $\lambda \not= 0$ are all smooth and isomorphic, whereas the zero fiber $M_n(0)$ is very singular.

\vskip 3mm

It is a consequence of deep results on motivic nearby cycles, due to Denef and Loeser (see \footnote{K. Behrend, J. Bryan and B. Szendroi, Motivic degree zero Donaldson-Thomas invariants, Inv. Math. 192, 111-160 (2013)}), that the {\em virtual} motive of the degeneracy locus is related to the difference of the motives of smooth and singular fibers
\[
[ \mathbf{rep}_n R_W ]_{virt} = \mathbb{L}^{- \frac{mn^2}{2}} ( [ M_n(0) ] - [M_n(1) ])  \quad \text{in} \quad \mathbf{Mot}_{\C}^{\mu_{\infty}}[ \mathbb{L}^{-\frac{1}{2}} ] \]
It takes a few words  to make sense of this equality.

\vskip 3mm

The ring of equivariant mtives $\mathbf{Mot}^{\mu_{\infty}}_{\C}$ is the ring of motives of pairs $(X,\mu_k)$ where $X$ is a reduced variety with an action of the cyclic group $\mu_d$ of $d$-th roots of unity, for some $d$. The scissor- and factoring-relations only hold in situations compatible with the group action. Further, we have that $[(\Af^n,\mu_d)]=\mathbb{L}^n$ whenever the action of $\mu_d$ on $\Af^n$ is linear. In the above equality we mean with $[M_n(1)]$ really $[(M_n(1),\mu_d)]$ where $\mu_d$ acts by multiplying each of the variables $X_i$. On the other hand, the action of $\mu_{\infty}$ on $M_n(0)$ is trivial.

\vskip 3mm

The square root of the Lefschetz motive $\mathbb{L}^{\frac{1}{2}}$ is the (equivariant) motive $1 - [ \mu_2 ]$, that is, the difference of the motive of one point with trivial action by two points which are interchanged under the $\mu_2$-action. To understand this we need a general trick on {\em separating variables}.

Consider two superpotentials $W \in \C \langle X_1,\hdots,X_m \rangle_d$ and $V \in \C \langle Y_1,\hdots,Y_k \rangle_d$ and denote the fibers of the corresponding  Chern-Simons functionals
\[
M_n(\lambda) = Tr(W)^{-1}(\lambda) \subset M_n(\C)^{\oplus m}  \quad \text{and} \quad N_n(\lambda) = Tr(V)^{-1}(\lambda) \subset M_n(\C)^{\oplus k} \]
and the fiber of the sum superpotential
\[
S_n(\lambda) = Tr(W+V)^{-1}(\lambda) \subset M_n(\C)^{\oplus m+k} \]
then we have the identity of equivariant motives
\[
[S_n(0)]-[S_n(1)] = ([M_n(0)]-[M_n(1)])([N_n(0)]-[N_n(1)]) \]
Indeed, we clearly have the 'formal' identity
\[
[S_n(\lambda)] = \sum_{\mu \in \C} [M_n(\mu)][N_n(\lambda-\mu)] \]
which can be made into a proper identity using the fact that $[M_n(\lambda)]=[M_n(1)]$ and $[N_n(\lambda)]=[N_n(1)]$ for all $\lambda \not= 0$. Then the above means
\[
\begin{cases}
[S_n(1)] = [M_n(0)][N_n(1)]+[M_n(1)][N_n(0)]+(\mathbb{L}-2)[M_n(1)][N_n(1)] \\
[S_n(0)] = [M_n(0)][N_n(0)] + (\mathbb{L}-1)[M_n(1)][N_n(1)]
\end{cases} 
\]
from which the identity of equivariant motives follows.

If we apply this to $W=X^2$ and $V=Y^2$, then we have for $n=1$ that
\[
[M_1(0)]-[M_1(1)]=[N_1(0)]-[N_1(1)]=1- [ \mu_2 ] \]
whereas for the sum potential $S=X^2+Y^2$ we easily compute
\[
[S_1(0)]=2 \mathbb{L}-1 \quad \text{and} \quad [S_1(1)]=\mathbb{L}-1 \]
and therefore we indeed have
\[
\mathbb{L}=[S_1(0)]-[S_1(1)]=([M_1(0)]-[M_1(1)])([N_1(0)]-[N_1(1)])=(1 - [ \mu_2])^2 \]

\vskip 3mm

Virtual motives are a bit harder to define properly (see \footnote{K. Behrend, Donaldson-Thomas type invariants via microlocal geometry, Ann. of Math. 170, 1307-1338 (2009)} for full details). If $X$ is a smooth variety then its virtual motive is a scaled version of the ordinary motive
\[
[X]_{virt} = \mathbb{L}^{-\frac{dim(X)}{2}} [X] \]
and in general the virtual motive depends on the singularities of the variety as well as on its embedding in $Y$ if it is the degeneracy locus of a functional $f : Y \rTo \C$.

\vskip 3mm

To appreciate the importance of the $\mu_d$-action, consider the superpotential  $W=XY-Z^2$.  Then $M_1(1) \subset \C^3$ is an affine smooth quadric with motive $\mathbb{L}^2+\mathbb{L}$ whereas $M_1(0) \subset \C^3$ is a singular quadric with motive $\mathbb{L}^2$. The degeneracy locus is a single point, the top of the cone in the zero-fiber. However, we get
\[
1 = [ \bullet]_{virt}  \overset{?}{=} \mathbb{L}^{-\frac{3}{2}} (\mathbb{L}^2-(\mathbb{L}^2 + \mathbb{L} )) = -\mathbb{L}^{-\frac{1}{2}} \]
This is caused by the fact that we didn't use the $\mu_2$-action on $M_1(1)$. If we take this action on $\mathbb{V}(xy-z^2-1)$ into account we can on the open piece where $x \not= 0$ again eliminate $y$ to get a contribution $\mathbb{L}(\mathbb{L}-1)$. But, on $x=0$ we have $y$ a free variable with $z=\pm i$ where the two values are interchanged under the $\mu_2$-action, so this piece contributes $\mathbb{L} [ \mu_2 ]$. In total we get
\[
[(M_1(1),\mu_2)] = \mathbb{L}^2-(1-[\mu_2])\mathbb{L} = \mathbb{L}^2-\mathbb{L}^{\frac{3}{2}} \]
and plugging this info in, we get a genuine identity.

\section{Brauer-Severis}

The fibers of the Chern-Simons functional can be viewed as (trace preserving) representation varieties, for we have
\[
M_n(\lambda) = \mathbf{trep}_n \mathbb{T}_n(\lambda) \quad \text{with} \quad \mathbb{T}_n(\lambda) =\frac{\mathbb{T}_{m,n}}{(Tr(W)-\lambda)} \]
with $\mathbb{T}_{m,n}$ the trace algebra of $m$ generic $n \times n$ matrices (that is, we adjoin to the ring generated by the generic matrices all traces of cyclic words in the generic matrices) and where we consider the trace preserving representations, that is, if $X_i$ is mapped to the matrix $A_i$, then its trace $Tr(X_i)$ must be mapped to $Tr(A_i)$.

\vskip 3mm

We can apply everything we did before, replacing $\mathbf{rep}_n R $ by $\mathbf{trep}_n R$, to the ring $R = \mathbb{T}_n(\lambda)$ in order to get
\[
(\mathbb{L}^n - 1) [ M_n(\lambda) ] = \sum_{k=1}^n [ S_{n,k}(\lambda) ] \quad \text{with} \quad S_{n,k}(\lambda) = \mathbf{S}_{n,k}(\mathbb{T}_n(\lambda)) \]
with
$[ S_{n,n}(\lambda) ] = [ BS_n(\lambda) ] [ GL_n] $ where $BS_n(\lambda) = \mathbf{BS}_n(\mathbb{T}_n(\lambda))$.
This time, the fibers of the map
\[
\psi : S_{n,k}(\lambda) \rOnto Gr(k,n) \]
over a $k$-dimensional subspace $V$ consists as before of points
\[
(\phi,v) = (\begin{bmatrix} \phi_1 & e \\ 0 & \phi_2 \end{bmatrix},\begin{bmatrix} w \\ 0 \end{bmatrix}) \]
But this time as the only relation is $Tr(\phi(W))=\lambda$, we have with $Tr(\phi_1(W)) = \mu$ that $Tr(\phi_2(W))$ must be equal to $\lambda - \mu$ and the extension $e$ can be arbitrary, giving us the formal  equality of motives
\[
[S_{n,k}(\lambda)] = \mathbb{L}^{m k (n-k)} [ Gr(k,n) ] \sum_{\mu \in \C} [ S_{k,k}(\mu) ] [ M_{n-k}(\lambda - \mu) ] \]
and as before we can convert this to a genuine identity. This then leads to  (see \footnote{Lieven Le Bruyn, Brauer-Severi motives and Donaldson-Thomas invariants of quantized threefolds . Journal of Noncommutative Geometry, 12, 671-692 (2018)}) 

\begin{theorem} The virtual motive $[\mathbf{rep}_n R_W]_{virt}$ can be computed inductively using the identity
\[
(\mathbb{L}^n - 1)([M_n(0)]-[M_n(1)]) = [GL_n]([BS_n(0)]-[BS_n(1)]) + \]
\[
 \sum_{k=1}^{n-1} \mathbb{L}^{mk(n-k)} [Gr(k,n)] [GL_k] ([BS_k(0)]-[BS_k(1)])([M_{n-k}(0)]-[M_{n-k}(1)]) \]
 That is, $[ \mathbf{rep}_n R]_{virt}$ can be computed from $[BS_k(0)]-[BS_k(1)]$ for all $1 \leq k \leq n$.
\end{theorem}

At first sight it might seem that computing $[M_n(\lambda)]$ is a lot easier than $[BS_n(\lambda)]$ as $M_n(\lambda)$ is a hyperplane in affine space
\[
M_n(\lambda) = \mathbb{V}(Tr(W)-\lambda) \subset \mathbf{rep}_n \C \langle X_1,\hdots,X_m \rangle = M_n(\C)^{\oplus m} \]
whereas $BS_n(\lambda)$ is a hyperplane in the generic Brauer-Severi variety
\[
BS_n(\lambda) = \mathbb{V}(Tr(W)-\lambda) \subset \mathbf{BS}_n(\C \langle X_1,\hdots,X_m \rangle) \]
Fortunately, Markus Reineke\footnote{M. Reineke, {\em Cohomology of non-commutative Hilbert schemes}, Alg. Rep. Thy. 8 (2005) 541-561} proved that $\mathbf{BS}_n(\C \langle X_1,\hdots,X_m \rangle)$ has a concrete cellular decomposition, with cells corresponding to sub-trees $\tau$ of the free $m$-ary tree consisting of $n$ nodes
\[
\mathbf{BS}_n(\C \langle X_1,\hdots,X_m \rangle) = \sqcup_{\tau} \mathbb{A}^{d(\tau)} \]
of which the dimensions $d(\tau)$ can be computed explicitly in terms of a right ordering of monomials in the $X_i$ corresponding to the nodes of the extended tree $\tilde{\tau}$ where we add leaves to all nodes of $\tau$.

\vskip 3mm

For example, $\mathbf{BS}_2(\C \langle X,Y \rangle) = \mathbb{A}^6 \sqcup \mathbb{A}^5$ where the two cells correspond to the sub-trees $\tau$ consisting of two nodes (solid edges) with the extended trees (dashed nodes)
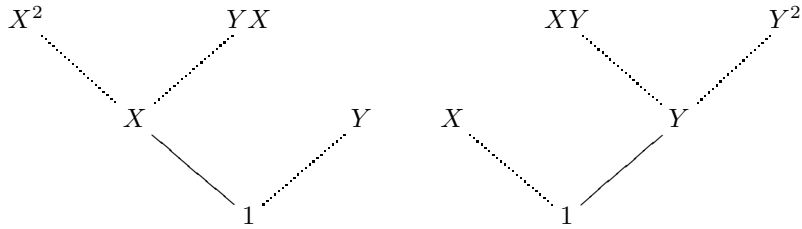
\[
\xymatrix{X^2 \ar@{.}[rd] & & YX \ar@{.}[ld] & \\
& X \ar@{-}[rd] & & Y \ar@{.}[ld] \\
& & 1 &} \qquad
\xymatrix{& XY \ar@{.}[rd] & & Y^2 \ar@{.}[ld] \\ X \ar@{.}[rd] & & Y \ar@{-}[ld] & \\ & 1 & &} \]
The monomials in these trees are ordered via
\[
\boxed{1} < \boxed{X} < X^2 < YX < Y \quad \text{and} \quad \boxed{1} < X < \boxed{Y} < XY < Y^2 \]
with the boxed terms the nodes of the tree. The dimension of the cell is then the sum over the extended leaves of the number of boxed terms which are smaller, that is, in our examples $2+2+2=6$ resp. $1+2+2=5$. This can then be used to give an explicit parametrization of the cells. Here, the two cells consist of triples $(X,Y,v)$ with
\[
( \begin{bmatrix} 0 & b \\ 1 & d \end{bmatrix},\begin{bmatrix} e & f \\ g& h \end{bmatrix},\begin{bmatrix} 1 \\ 0 \end{bmatrix}) \quad \text{resp.} \quad ( \begin{bmatrix} a & b \\ 0 & d \end{bmatrix},\begin{bmatrix} 0 & f \\ 1 & h \end{bmatrix},\begin{bmatrix} 1 \\ 0 \end{bmatrix})  \]
which makes it easy to compute the motive of the hypersurface $\mathbb{V}(Tr(W)-\lambda)$ in each cell by elimination of variables (we loose $n^2-n$ variables in going from $M_n(\lambda)$ to $BS_n(\lambda)$).

\vskip 3mm

As an example, consider the superpotential $W=X^2Y+YX^2$ with corresponding Jacobian algebra
\[
R_W = \frac{\C \langle X,Y \rangle}{(XY+YX,X^2)} \]
For starters, we have
\[
[M_1(0)]-[M_1(1)] = [BS_1(0)]-[BS_1(1)] = (2 \mathbb{L}-1) - (\mathbb{L}-1) = \mathbb{L} \]
giving $[\mathbf{rep}_1 R_W]_{virt} = 1$. To calculate $[\mathbf{rep}_2 R_W]_{virt}$ we need $[BS_2(0)]-[BS_2(1)]$ which we get from adding the contributions of each of the two cells in the generic Brauer-Severi variety.
The fiber $Tr(W)^{-1}(\lambda)$ in the first cell is the hyperplane
\[
H(\lambda) = \mathbb{V}(2 d^2h+2 bh + 2 bdg + 2 df + 2 be - \lambda) \subset \Af^6 \]
To compute $[H(0)]-[H(1)]$ we first consider the open piece where $d \not= 0$ on which we can eliminate $f$ independent of $\lambda$, so this piece does not contribute. If $d=0$ we have $2 bh + 2 be = \lambda$ with $f$ and $g$ free variables. On $b \not= 0$ we can eliminate $e$ independent of $\lambda$ so this does not contribute, and on $b=0$ we only get a contribution to $[H(0)]$ with $e,f,g$ and $h$ free variables, so $[H(0)]-[H(1)]=\mathbb{L}^4$.
The fiber $Tr(W)^{-1}(\lambda)$ in the second cell is the hyperplane
\[
H'(\lambda) = \mathbb{V}(2 d^2h + 2 bd + 2 ab - \lambda) \subset \Af^5 \]
(here, $f$ is a free variable).  On $b \not= 0$ we can eliminate $a$ independent of $\lambda$, so this does not contribute to $[H'(0)]-[H'(1)]$. If $b=0$ we have $d^2h = \lambda$ (with $a$ and $f$ free variables) giving $[H'(0)]-[H'(1)]=\mathbb{L}^3$. That is
\[
[BS_2(0)]-[BS_2(1)]=\mathbb{L}^4+\mathbb{L}^3 \quad \text{giving} \quad [\mathbf{rep}_2 R_W]_{virt} = \mathbb{L}^2 \]

\section{Another example}

As a (belated) answer to a question, we will compute the virtual representation motives of a speciafic {\em contraction algebra}, which in general are $2$-generated superpotential algebras $R_W$ assigned to divisorial contractions to curves in $3$-folds. 

Michael Wemyss\footnote{Will Donovan and Michael Wemyss, {\em Noncommutative enhancements of contractions} (2016), {\tt arXiv:1612.01687}} tells me there are reasons to conjecture that the virtual representation motives $[\mathbf{rep}_n R_W]_{virt}$ of such algebras are fully determined by those in small dimensions $n$. For more details see his paper and references contained in it.

\vskip 3mm

Consider the superpotential $W=X^3+Y^3$ with corresponding Jacobian algebra
\[
R_W = \frac{\C \langle X,Y \rangle}{(X^2,Y^2)} \]

\vskip 3mm

By separation of variables, we have for all $n$ that
\[
[\mathbf{rep}_n R_W]_{virt} = ( \mathbb{L}^{-\frac{n^2}{2}} ([M_n(0)]-[M_n(1)]) )^2 \]
with $M_n(\lambda) = \{ A \in M_n(\C)~|~Tr(A^3)=\lambda \}$. In this $m=1$ case, Reineke's decomposition result concerns a single tree
\[
\xymatrix{1 \ar@{-}[r] & A \ar@{-}[r] & A^2 \ar@{-}[r] & \hdots  \ar@{-}[r] & A^{n-1} \ar@{.}[r] & A^n} \]
which gives us that
\[
\mathbf{BS}_n(\C \langle A \rangle) = \mathbb{A}^n \]
with parametrization
\[
(A,v) = ( \begin{bmatrix} 0 & 0 & \hdots & 0 & a_1 \\
1 & 0 & \hdots & 0 & a_2 \\
0 & 1 & \hdots & 0 & a_3 \\
\vdots & \vdots & \ddots & \vdots & \vdots \\
0 & 0 & \hdots & 1 & a_n \end{bmatrix}, \begin{bmatrix} 1 \\ 0 \\ 0 \\ \vdots \\ 0 \end{bmatrix} ) \]
But then, for $n \geq 3$ we have that
\[
BS_n(\lambda) = \{ A~|~Tr(A^3)=\lambda \} = \mathbb{V}(a_n^3 + 3 a_{n-1}a_n + 3 a_{n-2}-\lambda) \subset \mathbb{A}^n \]
from which we can eliminate $a_{n-2}$ independent of $\lambda$, whence
\[
[BS_n(0)]-[BS_n(1)] = 0 \quad \text{for all $n \geq 3$} \]
and therefore $[\mathbf{rep}_n R_W]_{virt}$ for all $n \geq 3$ is fully determined by $[\mathbf{rep}_1 R_W]_{virt}$ and $[\mathbf{rep}_2 R_W]_{virt}$. More explicitly, we have
\[
[M_1(0)]-[M_1(1)]=[BS_1(0)]-[BS_1(1)]=1 - [ \mu_3 ] \]
and therefore
\[
[ \mathbf{rep}_1 R_W ]_{virt} = \mathbb{L}^{-1} (1 - [\mu_3])^2 \]
For $n=2$ we have
\[
BS_2(\lambda) = \mathbb{V}(a_2^3 + 3 a_1 a_2-\lambda) \subset \mathbb{A}^2 \]
On $a_2 \not= 0$ we can eliminate $a_1$ independent of  $\lambda$, so this piece does not contribute. On $a_2=0$ we have that $a_1$ is a free variable, but only when $\lambda=0$, so we get
\[
[BS_2(0)]-[BS_2(1)]= \mathbb{L} \quad \text{and} \quad [M_2(0)]-[M_2(1)]= \mathbb{L}^2 (\mathbb{L}-1) + \mathbb{L} (1 - [ \mu_3 ])^2 \]
from which we obtain that
 \[
 [\mathbf{rep}_2 R_W]_{virt} =( (\mathbb{L}-1) + \mathbb{L}^{-1} (1 - [ \mu_3 ])^2)^2 \]
 For $n \geq 3$ we have the recurrence relation $[M_n(0)]-[M_n(1)]=$
 \[
 \mathbb{L}^{n-1}(1 - [ \mu_3])([M_{n-1}(0)]-[M_{n-1}(1)] + \mathbb{L}^{2n-2}(\mathbb{L}^{n-1}-1)([M_{n-2}(0)] - [ M_{n-2}(1)]) \]
 
 \vskip 3mm
 
 Similarly, we have for the superpotential $W = X^d+Y^d$ that $[BS_n(0)]-[BS_n(1)] = 0$ for all $n \geq d$ and therefore that the virtual motives $[\mathbf{rep}_n R_W]_{virt}$ can be computed from those with $n \leq d-1$.
 
 \vskip 5mm
 
 \noindent
 {\bf Acknowledgements.} This is the write up of a talk given in the Glasgow algebra seminar. I thank Theo Raedschelders and Michael Wemyss for an enjoyable stay and stimulating discussions (leading to the final example).

\end{document}